\documentclass[a4paper]{article}
\usepackage{amsfonts}
\usepackage{amsmath}
\usepackage{amssymb}
\begin{document}
\newtheorem{theorem}{Theorem}[section]
\newtheorem{lemma}[theorem]{Lemma}
\newtheorem{corollary}[theorem]{Corollary}
\newtheorem{conjecture}[theorem]{Conjecture}
\newtheorem{remark}[theorem]{Remark}
\newtheorem{definition}[theorem]{Definition}
\newtheorem{problem}[theorem]{Problem}
\newtheorem{example}[theorem]{Example}
\newtheorem{proposition}[theorem]{Proposition}
\title{{\bf Canonical volume forms on compact K\"{a}hler manifolds}}
\date{June 29, 2007}
\author{Hajime TSUJI}
\maketitle
\begin{abstract}
\noindent We construct a canonical singular hermitian metric with semipositive 
curvature current on the canonical line bundle of a compact K\"{a}hler 
manifold with pseudoeffective canonical bundle.
The method of the construction is a modification of the one in  \cite{super}.  
\\  MSC: 14J15,14J40, 32J18
\end{abstract}
\section{Introduction}

In \cite{super}, I have constructed the canonical singular hermitian metric
on relative canonical bundles of projective families whose general fiber has 
pseudoeffective canonical bundles. 

Although the construction in \cite{super} works only for projective families, 
it is very likely that the similar result holds also for K\"{a}hler families. 
The purpose of this short article is to show that this is the case up to 
certain  extent, i.e., the construction of the metric also works for the K\"{a}hler case, but we need more work to show the semipositivity of the curvature of the 
singular hermitian metrics on the families.  The reason is that we cannot use the branched coverings in the case of K\"{a}hler families. 
To overcome this difficulty, it seems to 
be necessary to consider the variation of plurisubharmonic objects instead 
of holomorphic objects.  
This ``horizontal semipositivity'' of the curvature is very important 
in many applications such as the deformation invariance of plurigenera in 
the case of K\"{a}hler families. 

Let $X$ be a compact K\"{a}hler manifold of dimension $n$  and let $K_{X}$ be the canonical bundle of $X$.   The purpose of this short note is to consruct a canonical singular 
hermitian metric on $K_{X}$ when $K_{X}$ is pseudoeffective.

\begin{theorem}\label{main}
Let $X$ be a compact K\"{a}hler manifold with pseudoeffective $K_{X}$. 
Then there exists a singular hermitian metric $\tilde{h}_{can}$ on $K_{X}$ such that 
\begin{enumerate}
\item $\tilde{h}_{can}$ is uniquely determined by $X$. 
\item $\tilde{h}_{can}$ is an AZD of $K_{X}$, i.e., $\tilde{h}_{can}$ 
has the following properties : 
\begin{enumerate}
\item The curvature current $\Theta_{\tilde{h}_{can}}$ of $\tilde{h}_{can}$ is semipositive in the sense of current on $X$. 
\item For every $m\geqq 0$, 
\[
H^{0}(X,{\cal O}_{X}(mK_{X})\otimes{\cal I}(\tilde{h}_{can}^{m}))
\simeq H^{0}(X,{\cal O}_{X}(mK_{X}))
\]
holds for every $m\geqq 0$. 
\end{enumerate}
\end{enumerate}
$\square$
\end{theorem}
\begin{remark}
If $X$ is projective, $\tilde{h}_{can}$ seems to be closely related 
to the supercanonical AZD $\hat{h}_{can}$ in \cite{super} (see Remark \ref{relation} below). $\square$ 
\end{remark}

We  note that the bounded semipositive $(n,n)$ form 
$\tilde{h}_{can}^{-1}$ is an invariant volume form on $X$, i.e., the 
automorphism group of $X$ preserves $\tilde{h}_{can}^{-1}$. 
 
\section{Construction of $\tilde{h}_{can}$}
In this section, we shall prove Theorem \ref{main}.

\subsection{The construction}

The construction below is modeled after the one in 
\cite{super}.  

Let $X$ be a compact K\"{a}hler manifold. 
Let $\omega_{0}$ be any K\"{a}hler form on $X$.  
For every $\varepsilon \in (0,1]$, we set 
\[
h_{\varepsilon}:= \inf \{ h \mid \mbox{ $h$ is a singular hermitian metric 
on $K_{X}$}, \Theta_{h} + \varepsilon\cdot\omega_{0} \geqq 0, \int_{X}h^{-1} = 1\},
\]
where the infimum means the pointwise infimum. 
And set 
\[
dV_{\varepsilon}:= h_{\varepsilon}^{-1}
\]
And we set 
\[
\tilde{h}_{can} : = \mbox{the lower envelope of}\,\,\,\,\liminf_{\varepsilon\downarrow 0} h_{\varepsilon}
\]
and define 
\[
d\tilde{V}_{can}:= \tilde{h}_{can}^{-1} 
= \mbox{the upper envelope of}\,\,\,\,\limsup_{\varepsilon\downarrow 0}
dV_{\varepsilon}.
\]
Then since the upper limit of a sequence of plurisubharmonic function 
locally uniformly bounded from above is again plurisubharmonic, if we take 
the upper envelope of the limit (cf. (\cite[p.26, Theorem 5]{l}), 
we see that 
$\tilde{h}_{can}$ has semipositive curvature in the sense of current. 
\subsection{Upper estimate of $dV_{\varepsilon}$}
Let $d\mu$ be a $C^{\infty}$ volume form on $X$. 
\begin{lemma}\label{upper}
There exists a positive constant $C$ independent of $\varepsilon \in (0,1]$
such that 
\[
dV_{\varepsilon} \leqq C\cdot d\mu
\]
holds on $X$. $\square$
\end{lemma}
{\bf Proof of Lemma \ref{upper}}. 
Let $x$ be a point on $X$. 
Let $(U,z_{1},\cdots ,z_{n})$ be a local coordinate neighbourhood such such that\begin{enumerate}
\item $z_{1}(x)= \cdots = z_{n}(x) = 0$.
\item $z_{1},\cdots ,z_{n}$ are holomorphic on a neighbourhood of the closure of $U$.
\item $U$ is biholomorphic to the unit open polydisk $\Delta^{n}$ in $\mathbb{C}^{n}$ with center $O$ via the coordinate $(z_{1},\cdots ,z_{n})$. 
\end{enumerate} 
By the $\partial\bar{\partial}$ Poincar\'{e} lemma, there exists a 
positive $C^{\infty}$ function $\varphi_{U}$ such that 
\[
\omega_{0}\mid_{U}= \sqrt{-1}\partial\bar{\partial}\log \varphi_{U}
\]
holds.  By the construction, we may and do assume that $\varphi_{U}$ can be 
taken so that $\varphi$ is bounded on a neighbourhood of the closure of 
$U$.  

Let $dV$ be a bounded uppersemicontinuous semipositive $(n,n)$ form such that 
\[
\sqrt{-1}\partial\bar{\partial}\log dV + \varepsilon\omega_{0}
\geqq 0
\]
and 
\begin{equation}\label{normal}
\int_{X}dV = 1
\end{equation}
hold.
We define the function $a_{U}$ by 
\[
dV\mid U = a_{U}\cdot (\sqrt{-1})^{n}\prod_{i=1}^{n}dz_{i}\wedge d\bar{z}_{i}
\] 
holds. 
Then since $\varepsilon \in (0,1]$, we have that 
\begin{equation}\label{psh}
\sqrt{-1}\partial\bar{\partial}\log (a_{U}\cdot e^{-\varphi_{U}}) \geqq 0
\end{equation}
holds in the sense of current.
And by (\ref{normal})
\begin{equation}\label{int}
\int_{U}e^{-\varphi_{U}}dV \leqq \sup_{U}e^{-\varphi_{U}} 
\end{equation}
holds. 
Hence combining (\ref{psh}) and (\ref{int}), by the submeanvalue property of 
plurisubharmonic functions, we see that 
\[
dV(x) \leqq (\sup_{U}e^{-\varphi_{U}})\cdot e^{\varphi_{U}(x)}
(\frac{\sqrt{-1}}{2\pi})^{n}\cdot\prod_{i=1}^{n}dz_{i}\wedge d\bar{z}_{i}
\]
holds. 
Hence since $X$ is compact, moving $x$, we obtain Lemma \ref{upper}. $\square$

\subsection{Lower estimate of $dV_{\varepsilon}$}
The lower estimate is also easy. 
Let $h$ be any singular hermitian metric on $K_{X}$ with semipositive 
curvature current.  We note that $h^{-1}$ is a bounded semipositive $(n,n)$ from on $X$. 
Then by the definition of $d\tilde{V}_{\varepsilon}$, we see that 
\[
dV_{\varepsilon} \geqq (\int_{X}h^{-1})^{-1}\cdot h^{-1}
\]
holds.   
Hence letting $\varepsilon$ tend to $0$, we have the following lemma.
\begin{lemma}\label{lower} Let $h$ be any singular hermitian metric on $K_{X}$ 
with semipositive curvature current. 
Then we have that 
\[
d\tilde{V}_{can} \geqq (\int_{X}h^{-1})^{-1}\cdot h^{-1}
\]
holds. $\square$
\end{lemma}

Lemma \ref{lower} immediately implies the following lemma. 

\begin{lemma}\label{azd}
$\tilde{h}_{can}$ is an AZD of $K_{X}$, i.e.,
\begin{enumerate}
\item $\Theta_{\tilde{h}_{can}}$ is semipositive in the sense of current. 
\item For every $m\geqq 0$, 
\[
H^{0}(X,{\cal O}_{X}(mK_{X})\otimes {\cal I}(\tilde{h}_{can}^{m}))
\simeq H^{0}(X,{\cal O}_{X}(mK_{X}))
\]
holds. $\square$ 
\end{enumerate}
\end{lemma}
{\bf Proof of Lemma \ref{azd}}. 
Let $\sigma$ be a nonzero element of $H^{0}(X,{\cal O}_{X}(mK_{X}))$ 
for some $m\geqq 0$. 
Then we see that 
\[
\frac{1}{\mid\sigma\mid^{\frac{2}{m}}} := \frac{h_{0}}{(h_{0}^{m}(\sigma,\sigma))^{\frac{1}{m}}}
\]
is a singular hermitian metric on $K_{X}$ with semipositive curvature current,
where $h_{0}$ is any $C^{\infty}$ hermitian metric on $K_{X}$. 
Hence by Lemma \ref{lower}, we see that 
\[
\tilde{h}_{can} \leqq \mid\int_{X}\mid\sigma\mid^{\frac{2}{m}}\mid\cdot\frac{1}{\mid\sigma\mid^{\frac{2}{m}}}. 
\]
This means that 
\[
\tilde{h}_{can}^{m}(\sigma ,\sigma ) \leqq \mid\int_{X}\mid\sigma\mid^{\frac{2}{m}}\mid^{m}
\]
holds.   Hence 
\[
\sigma \in H^{0}(X,{\cal O}_{X}(mK_{X})\otimes {\cal I}(\tilde{h}_{can}^{m}))
\]
holds.   Since $\sigma$ is arbitrary,  we complete the proof of 
Lemma \ref{azd}. $\square$ 
\begin{remark}\label{relation}
The above proof of Lemma \ref{lower} and \cite[Lemma 2.5]{super} implies that if $X$ is projective,
 the supercanonical AZD $\hat{h}_{can}$ of $K_{X}$ constructed in \cite{super}
 satisifies the inequalities
\[
\tilde{h}_{can} \leqq \hat{h}_{can} \leqq (\int_{X}\tilde{h}_{can}^{-1})\cdot \tilde{h}_{can}
\]
on $X$. $\square$  
\end{remark}

\subsection{Independence from $\omega_{0}$}

\begin{lemma}\label{indep}
$d\tilde{V}_{can}$ is independent of the choice of $\omega_{0}$. 
$\square$
\end{lemma}
{\bf Proof of Lemma \ref{indep}}. 
Let $\omega_{0}^{\prime}$ be another K\"{a}hler form on $X$. 
Then since $X$ is compact, there exists a positive constant $C_{0} > 1$
such that 
\[
C_{0}^{-1}\cdot \omega_{0}^{\prime} \leqq \omega_{0} \leqq C_{0}\cdot\omega_{0}^{\prime}
\]
holds on $X$. 
For every $\varepsilon \in (0,1]$, we set 
\[
h^{\prime}_{\varepsilon}:= \inf \{ h \mid \mbox{ $h$ is a singular hermitian metric 
on $K_{X}$}, \Theta_{h} + \varepsilon\cdot\omega^{\prime}_{0} \geqq 0, \int_{X}h^{-1} = 1\},
\]
where the infimum means the pointwise infimum. 
And set 
\[
dV^{\prime}_{\varepsilon}:= (h_{\varepsilon}^{\prime})^{-1}. 
\]
Then by the definitions of $h_{\varepsilon}$ and $dV_{\varepsilon}$ 
we have the inequalities :
\[
h^{\prime}_{C_{0}\varepsilon}\leqq h_{\varepsilon} \leqq h_{C_{0}^{-1}\varepsilon}
\] 
and 
\[
dV^{\prime}_{C_{0}^{-1}\varepsilon}
\leqq dV_{\varepsilon} \leqq dV_{C_{0}\varepsilon}
\]
hold for every $\varepsilon \in (0,C_{0}^{-1}]$. 
Hence letting $\varepsilon$ tend to $0$, we see that 
$d\tilde{V}_{can}$ is independent of the choice of $\omega_{0}$. $\square$ 

\subsection{Completion of the proof of Theorem \ref{main}}
By  Lemmas \ref{upper} and \ref{lower} and the construction of 
$\tilde{h}_{can}$, we see that $\tilde{h}_{can}$ is a well defined 
singular hermitian metric on $K_{X}$. 
Also by Lemma \ref{azd}, we see that $\tilde{h}_{can}$ is an AZD of 
$K_{X}$.   By Lemma \ref{indep}, we see that $\tilde{h}_{can}$ is 
independent of the choice of $\omega_{0}$. 

Hence we complete the proof of Theorem \ref{main}. $\square$  
\subsection{A variant of the construction}
One may modify the construction above.  In fact one may define 
the singular hermitian metric without using the auxilary K\"{a}hler form 
$\omega_{0}$.
\[
\breve{h}_{can}:= \inf \{ h\mid \mbox{$h$ is a singular hermitian metric on $K_{X}$},\,\, \Theta_{h} \geqq 0, \int_{X}h^{-1} = 1\}. 
\] 
Then by the same proof as above, we may prove that $\breve{h}_{can}$ satisfies the 
same properties in Theorem \ref{main}. 
This construction is a little bit simpler.  
But it is not clear whether this construction has any advantage.

Author's address\\
Hajime Tsuji\\
Department of Mathematics\\
Sophia University\\
7-1 Kioicho, Chiyoda-ku 102-8554\\
Japan \\
e-mail address: tsuji@mm.sophia.ac.jp

\end{document}